\documentclass[12pt]{amsart}

 \usepackage{graphicx}
 \usepackage{epstopdf}
     \usepackage{amssymb}
  \usepackage{color} 

\long\def\symbolfootnote[#1]#2{\begingroup%
\def\thefootnote{\fnsymbol{footnote}}\footnote[#1]{#2}\endgroup}

    \newcommand\Def{\mathop{\rm Def}\nolimits}

\newcommand\bsm{ \begin{smallmatrix}}
\newcommand\bspm{ \left(\begin{smallmatrix}}
\newcommand\esm{\end{smallmatrix} }
\newcommand\espm{\end{smallmatrix} \right)}
\newcommand\bbm{\left[\begin{matrix}}
\newcommand\ebm{\end{matrix}\right]}
\newcommand\bcs{\begin{cases}}
\newcommand\ecs{\end{cases}}

 \newcommand{\lra}[1]{\left\langle#1\right\rangle}
\newcommand{\lrc}[1]{\left\{ #1\right\}}

\newcommand{\lrp}[1]{\left(#1\right)}

 \newcommand\wt[1]{\widetilde{#1}}

 \newcommand\ff{\mathfrak{f}}

 \newcommand\sidecong{\rotatebox{90}{$\cong$}}

  \newcommand\sideeq{\rotatebox{90}{$=$}}
   
    \newcommand\sidein{\rotatebox{90}{$\in$}}

\newcommand{\C}{{\mathbb{C}}}

\newcommand{\G}{\mathbb{G}}

\renewcommand{\P}{\mathbb{P}}

\newcommand{\R}{\mathbb{R}}

\newcommand{\Z}{\mathbb{Z}}

\newcommand{\cC}{{\mathscr{C}}}
\newcommand{\cD}{{\mathscr{D}}}

\newcommand{\cO}{{\mathscr{O}}}

\newcommand{\cX}{{\mathscr{X}}}

\newcommand\bpm{\begin{pmatrix}}
\newcommand\epm{\end{pmatrix}}

\newcommand\mhs{mixed Hodge structure}

\newcommand\Pic{{\mathop{\rm Pic}\nolimits}}

\newcommand{\Hg}{{\rm Hg}}

\newcommand\Rim{\mathop{\rm Im}\nolimits}

\newcommand\id{\mathop{\rm id}\nolimits}

\newcommand\Res{\mathop{\rm Res}\nolimits}

 \newcommand{\Alb}{{\rm Alb}}

 \newcommand{\Hom}{\mathop{\rm Hom}\nolimits}

\newcommand{\AJ}{{\rm AJ}}

 \renewcommand{\part}{\partial}
\newcommand{\la}{{\lambda}}

\newcommand{\Om}{{\Omega}}
\newcommand{\om}{{\omega}}

\newcommand\sig{\sigma}

\newcommand\bsl{\backslash}

\newcommand{\lab}{\label}

\newcommand\ccup{\mathop{\cup}\limits}

\newcommand\bmp[2]{\hbox{\begin{minipage}[c]{#1in}  #2 \end{minipage}}}

\newcounter{demo}[equation]

 \newtheoremstyle{mytheo}
  {3pt}
  {3pt}
  {\itshape}
  {}
  {\scshape}
  {:}
  {.5em}
  {}

\theoremstyle{mytheo}

\newtheorem*{thma}{Theorem A}
\newtheorem*{thmb}{Theorem B}

 \newtheoremstyle{subsect}
 {3pt}
  {3pt}
  {}
  {}
  {\it}
  {\upshape{:}}
  {.5em}
  {}
\theoremstyle{subsect}

\newtheoremstyle{note}
  {3pt}
  {3pt}
  {}
  {}
  {\bfseries}
  {:}
  {.5em}
  {}
\theoremstyle{note}

\newtheorem*{rem}{Remark}

\theoremstyle{remark}

\newcommand\xri[1]{\xrightarrow{#1}}

\newcommand\ol[1]{\overline{#1}}

 \usepackage[mathscr]{euscript}

 \renewcommand{\theequation}{\thesubsection.\arabic{equation}} 
 
\usepackage{amsmath}
  \usepackage[all]{xy}
  \usepackage{hyperref}
\usepackage{amsfonts}

   \usepackage{enumerate}
    \usepackage{amssymb}
 \usepackage[mathscr]{euscript}
 
\numberwithin{equation}{section}

\DeclareFontFamily{U}{mathx}{}
\DeclareFontShape{U}{mathx}{m}{n}{<-> mathx10}{}
\DeclareSymbolFont{mathx}{U}{mathx}{m}{n}
\DeclareMathAccent{\widecheck}{0}{mathx}{"71}
\newcommand\fZ{\mathfrak Z}
\newcommand\fX{\mathfrak X}

  \newcommand\bsp[1]{\begin{split} #1 \end{split}}
  
  \newcommand\beb{\begin{enumerate}[$\bullet$]}
  \newcommand\bebI{\begin{enumerate}[I.]}
  \newcommand\bebi{\begin{enumerate}[{\rm (i)}]}
  \newcommand\beba{\begin{enumerate}[{\rm (a)}]}
  \newcommand\bed{\begin{enumerate}[--]}
 \newcommand\eeb{\end{enumerate}}

    \newcommand\Gy{\mathop{\rm Gy}\nolimits}

   \newcommand\Ypd{Y^{(d)}}

  \theoremstyle{mytheo}

   \theoremstyle{note}
 \newtheorem{Examo}{Example}

 \renewcommand{\theequation}{\thesection.\arabic{equation}}

   \DeclareFontFamily{U}{mathx}{}
\DeclareFontShape{U}{mathx}{m}{n}{<-> mathx10}{}
\DeclareSymbolFont{mathx}{U}{mathx}{m}{n}
\DeclareMathAccent{\widecheck}{0}{mathx}{"71}

\newcommand\opsi{\psi} 
        
\begin{document}
      \title[Lagrangian interpretation of Abel-Jacobi mappings]{Lagrangian interpretation of Abel-Jacobi mappings associated to Fano threefolds}
 \author[]{Rodolfo Aguilar, Mark Green and  Phillip Griffiths}
 
\address{Department of Mathematics, University of Miami,\hfill\break\indent Coral Gables, FL 33146 \hfill\break\indent
{\it E-mail address}\/: {\rm aaguilar.rodolfo@gmail.com}
}
\address{
Department of Mathematics, University of California at\hfill\break\indent Los Angeles, Los Angeles, CA 90095\hfill\break\indent
 {\it E-mail address}\/: {\rm mlg@ipam.ucla.edu}}
 
 \address{Department of Mathematics, University of Miami,\hfill\break\indent Coral Gables, FL 33146, and \hfill\break\indent Institute for Advanced Study, Einstein Drive, Princeton, NJ 08540\hfill\break\indent
{\it E-mail address}\/: {\rm pg@ias.edu}}   
 
 \begin{abstract}
 Using the general framework due to Donagi-Markman \cite{DM} and Markushevich \cite{M} we shall derive an expression for the differential of Abel-Jacobi mappings on Fano threefolds.  This formula involves information normal to the Lagrangian submanifolds constructed in \cite{DM} and \cite{M}.  It may be applied to give new proofs of a number of classical results about these varieties.
 \end{abstract}
        \maketitle
       \renewcommand{\theequation}{\thesection.\arabic{equation}} 
       
   \tableofcontents
             
       \section{Introduction; statement  of Theorem A}
       In this note we will be interested in the differential of Abel-Jacobi maps
       \begin{equation}\lab{1.1}
       \AJ_X:Z\to J(X)\end{equation}
       where $X$ is a smooth Fano threefold, $Z$   the parameter space of a family of  degree $d$   curves $C_z\subset X$,  $z\in Z$,  and 
       $$J(X)=H^1(\Om^2_X)^\ast / H_3(X,\Z)\cong H^2(\Om^1_X)/H^3(X,\Z)$$ is the intermediate Jacobian of $X$. The objective is  to analyze the differential $d\AJ_X$ of \eqref{1.1} in terms of the intersections $C_z\cdot Y$ where $Y\in |-K_X|$ is a smooth anti-canonical divisor.  For simplicity of exposition we will assume these intersections to be transverse.  They then induce a morphism 
       \[
       f:Z\to Y^{(d)}\]
       where the $d^{\rm th}$ symmetric product is a hyperK\"ahler manifold with symplectic form $\opsi\in H^0(\Om^2_{Y^{(d)}})$.  It is known (\cite{IM}), and will be verified below, that $f(Z)\subset Y$ is a maximal dimension Lagrangian submanifold.  We will define maps
       \begin{equation}\lab{1.2}\bsp{
       \alpha: &T_z Z\to T\Ypd\big|_{C_z\cap Y},\quad \alpha=f_\ast,\\
       \beta:& H^1(\Om^2_X)\to T\Ypd\big|_{C_z\cap Y}/\Rim \alpha; }\end{equation}
      thus,  $\beta$  maps to the normal space $N_{f(Z)/Y^{(d)},f(z)}$.        \begin{thma}
      {\rm (i)} $f(Z)$ is a maximal dimensional Lagrangian submanifold of $\Ypd$. {\rm (ii)} For $\xi\in T_zZ$ and $\eta\in H^1(\Om^2_X)$
      \[
      \opsi(\alpha(\xi),\beta(\eta))=\lra{d\AJ_X(\xi),\eta} \]
      where $\lra{\enspace,\enspace}$ is the pairing $H^2(\Om^1_X)\otimes H^1(\Om^2_X)\to \C$.
        \end{thma}
      We note that the above holds for any smooth $Y\in |-K_X|$ for which the intersections $C_z\cdot Y$ are transverse.  We also will see that $\beta$ and $\psi$ are each defined up to a constant and these are linked, so that the formula in Theorem A makes sense.
 The result will enable us to analyze $d\AJ_X$ in terms of the symplectic geometry associated to the pair $(X,Y)$.

As an application of this result we will give a  simple proof of the tangent bundle theorem (TBT) for the Fano surface $Z$ of lines on a smooth cubic threefold.  The proof of this   central  result will have as one consequence  that the induced map $\Alb(Z)\to J(X)$ is an isomorphism.  
      
  Theorem B will deal with the case where we deform $X$ subject to the constraint that $X$ contain a  $Y\in |-K_X|$ that is fixed in moduli as we deform $X$.   As will be recalled below the deformation theory for this is unobstructed with Kuranishi space an open set $B\subset H^1(\Om^2_X)$.  The first order variation of $X$  appears in the definition of $\beta$ in \eqref{1.2}.
      
      Given the family of curves $C_z\subset X$, $z\in Z$, there is a corresponding deformation theory for  the family.  We assume that this theory is also unobstructed and, as will be explained below, this leads to the diagram    \begin{equation}\lab{1.4}\bsp{
 \xymatrix@C=.25pt{\fZ\ar[d]_{\pi_\fZ}& &\fX \ar[d]\\
      B_Z  &\subset &B}\lower17pt\hbox{$ ,\quad \cC\subset \fZ\times_{B_Z}\cX$}}\end{equation}
      where $B_Z$ is an open set in $\ker\{ H^1(\Om^2_X)\to H^1(N_{C/X})\}$.  We will assume that $\fZ$ and $B_Z$ are smooth and that $\fZ\to B_Z$ is a smooth fibration.  Corresponding to points $(b,z)\in\fZ$ where $b\in B_Z$ and $z\in Z_b=\pi^{-1}_\fZ(b)$ there is a curve $C_{b,z}\subset X_b$ and this defines $\cC\subset \fZ\times_{B_Z}\cX$ in \eqref{1.4}.
      
      The mapping of the curves $C_{b,z}$ to intermediate Jacobians is accomplished by using normal functions and Deligne cohomology \cite{DM}.  For the latter there is for each $X_b$ a Deligne cohomology group $H^2(\cD(X_b))$ to which the codimension-2 algebraic cycles on $X_b$ are mapped.  There is an exact sequence
      \begin{equation}\lab{1.5}
      0\to J(X_b)\to H^2(\cD(X_b))\to \Hg^2(X_b)\to 0\end{equation}
      where the last map sends an algebraic cycle to its fundamental class, and where the kernel of that map is the Abel-Jacobi map.  We denote by $\cD_\la(X_b)$ the subgroup of $H^2(\cD(X_b))$ that maps to the class $[C_{b,z}]\in H^4(X_b,\Z)$ of $C_{b,z}$ or to its translate under monodromy.  This is defined over   $B_Z\subset B$.  Letting $\cD_\la=\ccup_{b\in B_Z}\cD_\la(X_b)$ we define
 \begin{equation}\lab{1.6}
 F:\fZ\to Y^{(d)}\times \cD_\la\end{equation}
 by
 \[
 F(b,z) = (C_{b,z}\cap Y,\; [C_{b,z}]_{\cD}).\]
  
 By generalizing the constructions in \cite{DM} to the relative case of $(X,Y)$'s, it may be proved  that there is a holomorphic symplectic structure with symplectic form $\tau$ on $\cD_\la$ such that $\cD_\la\to B$ is a Lagrangian fibration (cf.\ \cite{M}).  We denote by
 \[
 \sig=\opsi+\tau\]
 the holomorphic symplectic form on $Y^{(d)}\times \cD_\la$.
      
      \begin{thmb}
      $F(\fZ)$ is a maximal dimension Lagrangian submanifold in $Y^{(d)}\times \cD_\la$.
      \end{thmb}
      
 What will be shown below  for part of the proof of Theorem B is
      \begin{equation}
      \lab{1.7} 
 \dim F(\fZ)=\frac12 \dim(Y^{(d)}\times \cD_\la) = d+h^{2,1}(X)\end{equation}
  The remainder of the argument may be done by an adaptation of the arguments in \cite{DM}. Both (i) in Theorem A and Theorem B are special cases of a   general result  whose proof serves to   illuminate the structures involved.  This will be carried out elsewhere.
     
     For simplicity of exposition we will assume that $H^0(T_X)=0$. Most of the Fano varieties that have non-trivial intermediate Jacobians will satisfy this condition.
     
     There is an extensive literature dealing with the topics in this note.  Beginning with \cite{DM} and continuing through the present time, Lagrangian fibrations have been a central subject in algebraic geometry and in its connection to physics.  Similarly, the study of Fanos, Fano pairs and the associated Abel-Jacobi mappings has been a topic of major and ongoing interest.   This is especially the case when singularities are allowed.  In this  note we  provide in the smooth case  one connection among the above areas.

      \section{Deformations of  anti-canonical pairs\\ and proof of {\rm (i)} in Theorem {\rm A}}

       In the proof of Theorem A we will use the deformations of the above structure (cf.\ \cite{B}).  The Kuranishi space of deformations of $X$ keeping $Y$ fixed is realized as an open set $B\subset H^{2,1}(X)$.     We have the mapping
      \begin{equation}\lab{new2.1}
      \xymatrix@R=1pt{\ff:\fZ\ar[r]&Y^{(d)}\\
      \sidein&\sidein\\
      (b,z)\ar[r]&C_{b,z}\cdot Y }\end{equation}
      which is the first factor in the mapping \eqref{1.6} defined above.
      
      We now turn to the proof of (i) in Theorem A.  There are several basic exact sheaf sequences and much of the argument here and later consists in using duality for the groups in their   cohomology sequences and in some cases non-trivially proving that maps that arise are the natural ones suggested by the notations.  The one exception to this is that we crucially use the non-degeneracy in both variables of the symplectically  induced bilinear form $\psi:T_{f(z)}f(Z)\otimes N_{f(Z)/\Ypd,f(z)} \to \C$.
      
      In what follows we fix a generator $\Psi\in H^0(-K_X)\cong H^0(\Om^3_X(Y))$, which then gives
      \begin{equation}\lab{2.1}
      T_X(-Y)\cong \Om^2_X.\end{equation}
      The first two exact sheaf sequences are
      \begin{equation}\lab{2.2}
      \bcs 0\to T_X(-\log Y)\to T_X\to N_{Y/X}\to 0\\
      0\to T_X(-Y)\to T_X (-\log Y) \to T_Y\to 0.\ecs\end{equation}
      Using that $-K_X$ is ample and that $Y$ is K3 surface,  the Akizuki-Nakano vanishing theorem and our assumption  $h^0(T_X)=0$ give for the cohomology sequences of \eqref{2.2}   the diagram
      \begin{equation}\lab{2.3} 
      \xymatrix{&&0\ar[d]&&&\\
     & &H^0(N_{Y/X})\ar[d]&&&\\
    0\ar[r]& H^1(\Om^2_X)\ar[r]& H^1(T_X(-\log Y))\ar[d]_{d\pi_X}\ar[r]^{\qquad d\pi_Y}& H^1(T_Y)\ar[r]&H^2(\Om^2_X)\ar[r]&0\\
    &&H^1(T_X)\ar[d]&&&\\
    &&0.&&&
    	}
	 \end{equation}  
	The mappings $d\pi_X, d\pi_Y$  are interpreted as mappings on tangent spaces in the diagram
	\begin{equation}\lab{2.4}
	 {
	\xymatrix@C=1pc{&\Def (X,Y)\ar[dr]^{\pi_Y}\ar[dl]_{\pi_X}&\\
	\Def(X)&&\Def(Y,R)}}\end{equation}
	where $R=\Rim\{\Pic(X)\to \Pic(Y)\}$ and $\Def(Y,R)$ is the deformation space of $Y$ preserving $R$.
	
	Letting $C$ denote a typical $C_z$, the second sheaf sequences are
	\begin{equation}\lab{2.5}\bsp{
	\lower.7in\hbox{$\left\{\begin{matrix}   \\[24pt] \\[34pt] \\[4pt] \end{matrix}\right. $}\xymatrix@R=.25pt@C=1pc{
	0\ar[r]& N_{C/X}\otimes K_X\ar[r]& N_{C/X}\ar[r]& N_{C\cap Y}\big|_Y \ar[r]& 0\\ 
	&&&\sideeq&\\ &&&T\Ypd\big|_{C\cap Y}&\\
	0\ar[r]& N_{C/X}\ar[r]& N_{C/\fX}\ar[r]& N_{X/\fX}\big|_C\ar[r]& 0\\ 
&&&\sidecong\\&&& H^{2,1}(X)\otimes \cO_C.} }\end{equation}
In the first sequence we are assuming that $C$ meets $Y$ transversely; this gives the vertical equality. In the second sequence are using the horizontal sequence in \eqref{2.3} to identify $N_{X/\fX}$.

\begin{proof}[Proof of {\rm (i)} in Theorem A]\hspace*{-5pt}\footnote{The argument  is essentially that given in \cite{IM}.}
  The cohomology sequences of  the first sequence in \eqref{2.5} and duality give
\[
T(\text{fibre $Z_f$ of }Z\xri{f} Y^{(d)})\cong H^1(N_{C/X})^\ast.\]
It follows that
\begin{align*}
\dim f(Z)&= h^0(N_{C/X})-h^1(N_{C/X})\\
&=\chi(N_{C/X})\\
&= \deg(-K_X\big|_C)=d.\end{align*}
The last step uses Riemann-Roch for $N_{C/X}\to C$, which turns out to have a special form due to $\dim X=3$ (cf.\ \cite{M}).

To show that $f(Z)$ is Lagrangian we  use $\wedge^2 T_zZ=\wedge^2 H^0(N_{C/X})$.  Locally we choose coordinates $u,v,w$ on $X$ so that around $p\in C\cap Y$ we have for $\xi,\eta\in   H^0(N_{C/X})$ 
\[\bcs
Y&\hspace*{-8pt}= \,\{w =0\},\\
\Psi&\hspace*{-8pt}= \,g(u,v,w)du\wedge dv\wedge dw / w,\\
\opsi&\hspace*{-8pt}=\,g(u,v,0)du\wedge dv\\
C&\hspace*{-8pt}=\,\{0,0,w)\}\\
\xi&\hspace*{-8pt}=\, \part/\part u,\; \eta=\part/\part v.\ecs\]
Then along $C$
\[
\xi\wedge \eta \rfloor \Psi \in H^0  (\Om^1_C(\log C\cap Y))\]
is well defined, and  for $\psi_f=f^\ast \psi$ we have from the above local expression
\[
\lra{\opsi_f,\xi\wedge\eta}=\sum_{p\in C\cap Y} \mathrm{Res}_p(\xi\wedge\eta\rfloor \Psi)=0.\qedhere\]
\end{proof}

\begin{rem} We remark that behind this argument is the general fact that in the exact cohomology sequence for the \mhs\ on $H^\ast(A\bsl B)$ for a smooth pair $A\subset B$ the composition $\Gy\circ\Res=0$ where $\Gy$ and $\Res$ are the Gysin and residue mappings.  This will be brought out in a second proof given below of (i) in Theorem A. 
\end{rem}
\section{Proof of {\rm (ii)} in Theorem {\rm A}}
There are two basic ideas behind the argument.  For the first we have for $\ff:\fZ\to \Ypd$ given by (\ref{new2.1})
\[
\ff^\ast\opsi\big|_Z=0.\]
Along $Z$, but \emph{not} restricted to the tangent spaces to $Z$, we have
\begin{equation}\lab{3.1}
\Om^2_\fZ\big|_Z\to \Om^2_Z\to 0 \end{equation}
and $\ff^\ast\opsi\big|_Z$ is in the kernel of this map.  Then
\begin{equation}\lab{3.2}
\ker\lrc{\Om^1_\fZ\to \Om^1_Z} \otimes \Om^1_Z\cong N^\ast_{Z/\fZ}\otimes \Om^1_Z\end{equation}
maps onto  the kernel in \eqref{3.1}.  Here $N^\ast_{Z/\fZ}$ is the conormal bundle to  $Z$ in $\fZ$.  Since the normal bundle to $Z$ in $\fZ$ is a sub-bundle of the trivial bundle with fibre $H^1(\Om^2_X)$ this brings $H^{2,1}(X)$ into the picture.

The second idea is that $\Psi\in H^0(\Om^3_X(\log Y))$ lifts to a section $\Psi_0\in H^0(\Om^3_{\fX/B_Z}(\log (Y\times B_Z)))$.  Here we are shrinking $B_Z$ to a neighborhood of the origin $b_0$.  If we could lift $ \Psi_0$ to a section in $H^0(\Om^3_\fX(\log (Y\times B_Z)))$, then the proof of (i) in Theorem A could be adapted to show that $\opsi$ vanishes on \eqref{3.2}, which is impossible since $f(Z)$ is maximal Lagrangian in $\Ypd$.  The obstruction to this lifting is the differential 
\[
T_{b_0} B_Z\to \Hom(H^0(\Om^3_X(\log Y)), H^1(\Om^2_X)) \]
of the period map in which once again $H^1(\Om^2_X)$ appears.  Here we are fixing $\Res\Psi_0=\psi\in H^0(\Om^2_Y)$; this implies that the differential of the period mapping maps
to 
$H^1(\Om^2_X)\subset H^1(\Om^2_X(\log Y))$.   Connecting these two ideas  is the central part of the proof.  

The sheaf $\Om^2_\fZ\big|_Z$ has a filtration with three graded pieces.  The argument below will deal with the most important middle one; the remaining one $\wedge^2 N^\ast_{Z/\fZ}$ may be treated directly and will also be analyzed  in the more general result to appear separately.

We first note that $\Psi_0\in H^0(\Om^3_{\fX/B_Z}(\log (Y\times B_Z)))$ may be lifted to a $C^\infty$ form
\[
\wt \Psi\in C^\infty(\Om^3_{\fX}(\log(Y\times B_Z))).\]
If $h=0$ is a local holomorphic equation of the divisor $Y\times B_Z$ in $\pi^{-1}_\fX(B_Z)$, then both $h \wt \Psi$ and $h\, d\wt \Psi$ are $C^\infty$.  We may take residues of forms in the $C^\infty$ log complex and the commutator $[\ol\part,\Res]=0$.

Next, using the present notation we will revisit the proof of (i) in Theorem A.  Setting $C=C_z$ and in the horizontal maps restricting along $C$  while in the vertical maps restricting to $X\subset \fX$ and $Z\subset \fZ$ we have
\[
\xymatrix{
\wt \Psi\in C^\infty(\Om^3_{\fX}(\log (Y\times B_Z)))\ar[d]\ar[r]& \Om^2_\fZ\otimes C^\infty (\log(\Om^1_C(C\cap Y)))\ar[d]\\
\Psi\in H^0(\Om^3_X(\log Y))\ar[r]& \Om^2_Z\otimes H^0(\Om^1_C(\log(C\cap Y))).}\]
Taking residues gives 
\[
\xymatrix{
\Om^2_\fZ\otimes C^\infty(\Om^1_C(\log (C\cap Y)))\ar[d]\ar[r]^{\quad\Res}&
 \Om^2_\fZ\otimes H^0(C\cap Y)(-1)\ar[r]^{\quad\Gy}\ar[d]& \ar[d]\Om^2_\fZ\otimes H^1(\Om^1_C)\\
\Om^2_Z\otimes H^0(\Om^1_C(\log(C\cap Y)))\ar[r]^{\quad\Res}& \Om^2_Z\otimes H^0(C\cap Y)(-1)\ar[r]^{\quad\Gy}& \Om^2_Z\otimes H^1(\Om^1_C)}\]
where $\Gy$ is the Gysin map and where the composition of the bottom arrows is zero.  This is another proof of  (i) in Theorem A.

Next we turn to the differential of the period mapping.  From the diagram
\[
0\to \Om^1_{B_Z}\otimes \Om^2_X (\log Y) \to \frac{\Om^3_\fX(\log(Y\times B_Z))\big|_X}{\Rim\{ \Om^2_{B_Z}\otimes \Om^1_\fX(\log (Y\times B_Z))\}}\to \Om^3_X(\log Y)\to 0\]
the connecting map in the exact cohomology sequence induces 
\[
\xymatrix@C=.1pc{
\Psi\ar[d]&\in&H^0(\Om^3_X(\log Y))\ar[d]\\
\nabla\Psi&\in& \Om^1_{B_Z} \otimes H^1(\Om^2_X(\log Y)).}\]
To compute this we lift $\Psi$ to $\wt\Psi$ as above and take $\ol\part\wt \Psi$. Because in the deformation $\fX\to B_Z$ the divisor $Y\times B_Z$ is fixed it follows that in the above we end up in $H^1(\Om^2_X)$.  Thus
\[
\ol\part\wt \Psi=i_\ast (\nabla \Psi),\]
where 
\[
i: \Om^1_{B_Z}\otimes \Om^2_X(\log Y) \to \Om^3_\fX(\log (Y\times B_Z))\big|_X.\]
Summarizing, the differential of the period mapping is
\begin{equation}\lab{3.3}\bsp{
\begin{picture}(250,120)
\put(0,100){$\Psi\in H^0(\Om^3_X(\log Y))$}
\put(5,90){\vector(0,-1){80}}
\put(0,0){$\nabla\Psi\in \Om^1_{B_Z}\otimes H^1(\Om^2_X).$}
\put(50,90){\vector(0,-1){20}}
\put(30,50){$\Om^1_{B_Z}\otimes H^1(\Om^2_X(\log Y))\xri{i_\ast}A^{0,1}(\Om^3_\fX(\log (Y\times B_Z))\big|_X)$}
\put(50,25){$\cup$}
\end{picture}}\end{equation}

We next do  for $\ol\part\wt\Psi$  the analogue of the above calculation  where we revisited the proof of (i) using $\Psi$.  For this we use the diagram
\begin{equation}\lab{3.4}\bsp{
\begin{picture}(320,120)
\put(0,100){$\wt\Psi\to \Res_{Y\times B_Z}\wt\Psi\in \Om^2_Y\to \Res_{Y\times B_Z}\wt\Psi\big|_C\in \Om^2_Y\big|_C\to \rho\lrp{\Res_Y\Psi} \in \Om^2_\fZ\otimes H^0(\cO_{C\cap Y})$}
\put(5,90){\vector(0,-1){20}}
\put(340,90){\vector(0,-1){20}}
\put(-20,55){$\ol\part\wt\Psi=i_\ast (\nabla\wt\Psi)\in \Om^1_{B_Z}\otimes H^{2,1}(X)$}
\put(280,55){$\ker(\Om^2_\fZ\to \Om^2_Z)\otimes H^1(\Om^1_C)$}
\multiput(40,45)(300,0){2}{\vector(0,-1){20}}
\put(20,5){$\Om^1_{B_Z}\otimes H^1(\Om^2_X)\big|_C$}
\put(290,5){$\Om^1_{B_Z}\otimes \Om^1_Z\otimes H^1(\Om^1_C)$}
\put(150,5){$\xri{\quad\id_{B_Z}\otimes \tau\quad}$}
\end{picture}}\end{equation}
where $\rho$ is the restriction mapping  $$\Om^3_\fX(\log (Y\times B_Z))\to \Om^2_\fZ\otimes \Om^1_C(\log(C\cap Y)),
$$
and $\tau$ is induced by considering a class in $H^1(\Om^2_X)\big|_C$ as a $(1,1)$ form in $A^1(\Om^1_C)$ whose values are a  conormal vector field  to $C$ in $X$.

We will  show below  that this diagram commutes.  Going by the upper path in it we obtain
\begin{equation}\lab{3.5}
\opsi(\alpha(\xi),\beta(\eta))\end{equation}
where $\xi\in H^0(N_{C/X})$, $\eta\in T_{b_0}B_Z\subset H^1(\Om^2_X)$.  Going by the lower path we have 
\begin{equation}\lab{3.6}
\lra{\nabla_\eta\Psi,\AJ_X(\xi)}\end{equation}
which is the desired result.

The commutativity of \eqref{3.4} uses the commutative diagram
\[
\xymatrix{ \Om^1_{B_Z}\otimes \Om^2_X(\log Y)\big|_C\ar[r]^{\id_{\Om^1_{B_Z}}\otimes \tau\qquad}\ar[d]^{i_\ast} & \Om^1_{B_Z}\otimes \Om^1_Z\otimes \Om^1_C(\log (C\cap Y))\ar[d]\\
\Om^3_\fX(\log(Y\times B_Z))\big|_C\ar[r]^\rho& \Om^2_\fZ\otimes \Om^1_C(\log (C\cap Y))}\]
which gives
\vspace*{3pt}

\centerline{\hfill$
\ol\part\rho (\wt\Psi\big|_C)=\rho (\ol\part\Psi\big|_C)=
\lrp{\id_{\Om^1_{B_Z}}\otimes \tau}(\nabla\Psi).\hfill\qed$}

\section{Lagrangian interpretation  of   Abel-Jacobi mappings; examples}
Together with duality the exact cohomology sequences of  \eqref{2.1} and \eqref{2.2}  give
\begin{align}\lab{4.1}\bsp{ 
0&\to H^0(N_{C/X})\to H^0(N_{C/\fX})\to H^{2,1}(X)\\&\xri{\delta}H^1(N_{C/X}) \to H^1(N_{C/\fX})
\to H^{2,1}(X)\otimes H^1(\cO_C)\to 0,}\end{align}
\begin{equation}\lab{4.2}\bsp{
\hspace*{-16pt}\xymatrix@R=.25pc@C=1pc{ 0\ar[r]&H^0(N_{C/X}\otimes K_X) \ar[r]&H^0(N_{C/X})\ar[r]^\alpha&
T\Ypd \big|_{C\cap Y}\ar[r]&H^1(N_{C/X}\otimes K_X)\ar[r]& H^1(N_{C/X})\ar[r]&0.\\
&\sidecong&&&\sidecong&&\\
&H^1(N_{C/X})^\ast&&& H^0(N_{C/X})^\ast&&}}\end{equation}
We will see below that $d\AJ_X\ne 0$ on the first term $H^1(N_{C/X})^\ast$  whereas obviously   $\alpha$ maps this space  to 0 in $T\Ypd\big|_{C\cap Y}$.
We also  recall that $\alpha=f_\ast$ in \eqref{4.2}.  The symplectic form $\opsi$ gives
\begin{equation}\lab{4.3}
H^1(\Om^2_X)\xri{\beta} T\Ypd\big|_{C\cap Y}\big/ \Rim\alpha.\end{equation}
Part (ii) in Theorem A gives for $\xi\in T_zZ\cong H^0(N_{C/X})$ and $\eta\in T_{b_0}B\cong H^1(\Om^2_X)$
\begin{equation}\lab{4.4}
\lra{d\AJ_X(\xi),\eta} = \lra{\alpha(\xi),\beta(\eta)}.\end{equation}
Finally the diagram
\begin{equation}\lab{4.5}\bsp{
\xymatrix@C=.2pc{\ker\alpha = H^0(N_{C/X}\otimes K_X)\ar[drr]_{d\AJ_X}&\cong& H^1(N_{C/X})^\ast\ar[d]^{\delta^\ast}\\
&&H^{1,2}(X)}}\end{equation}
commutes.

These sequences imply the interpretation
\begin{enumerate}[(1)]
\item $\xymatrix{\ker\alpha\ar[r]\ar[dr] &H^{1,2}(X)\ar[d]\\
& \Om^1_{B_Z};}$
 \item the failure of $C$ to freely deform with $X$ keeping $Y$ fixed, as measured by the codimension of $B_Z$ in $B$, is exactly compensated for by
\[
d\AJ_X\big|_{\ker\alpha}:H^1(N_{C/X})^\ast \to H^{1,2}(X);\]
in fact
\[
d\AJ_X\big|_{H^0(N_{C/X}\otimes K_X)} = \delta^\ast\]
where $\delta^\ast$ is given by \eqref{4.5}.\end{enumerate}

For a later application  we will use the identification
\begin{equation}\lab{4.6}
\ker(d\AJ^\ast_X) = \ker\beta \subseteq T_B\subseteq H^{2,1}(X).\end{equation}

\begin{proof}
From \eqref{4.1} and \eqref{4.2} we have
\[
\xymatrix{
0\ar[d]&\\
TB_Z\ar[d]\ar[r]^{\hspace*{-.85in}\beta}& \ker\{H^0(N_{C/X})^\ast \to H^1(N_{C/\fX})\}\ar[d]\\
H^{2,1}(X)\ar[dr]^\delta&H^0(N_{C/X})^\ast\ar[d]\\
&H^1(N_{C/X}).}\]
Using the duality given by $\psi$ we have
\begin{align*}
0\to \frac{H^0(N_{C/X})}{H^0(N_{C/X}(-Y))}\to T\Ypd\big|_{C\cap Y} &\to
\lrp{\frac{H^0(N_{C/X})}{H^0(N_{C/X}(-Y))}}^\ast\to 0,\\
\ker \lrp{ TB_Z\xri{\beta}\frac{TY^{(d)}\big|_{C\cap Y}}{\Rim\alpha}}& \cong \ker (d\AJ^\ast_X)\cap TB_Z\end{align*}
and
\[\xymatrix{
H^0(N_{C/X}(-Y))\ar[r]^{d\AJ_X}_{\begin{smallmatrix}\sideeq\\ \delta^\ast\end{smallmatrix}}&\lrp{\frac{H^{2,1}(X)}{TB_Z}}^\ast \to 0}\]
which gives \eqref{4.6}.
\end{proof}

\begin{Examo}
Let $X\subset \P V^\ast\cong \P^4 $ be a cubic threefold with Fano surface $Z\subset G(2,V^\ast) = \G(1,\P^4)$ where $G(2,V^\ast)$ is the Grassmannian of lines in $\P V^\ast$.  Denote by $C=C_z$ the line corresponding to $z\in Z$.

By \eqref{4.2}, the kernel of $\alpha:H^0(N_{C/X})\to TY^{(d)}\big|_{C\cap Y}$ is given by $H^1(N_{C/X})^\ast$.  But for the cubic 3-fold case $N_{C/X}\cong \cO_C\oplus \cO_C$ or $\cO_C(-1)\oplus \cO_C(1)$; therefore $H^1(N_{C/X})=0$ and $\alpha$ is injective.

To study the map
\[
\beta:H^1(\Om^2_X) \to TY^{(d)}_{C  \cap Y}/\Rim\alpha\cong N_{f(Z)/Y^{(d)},f(z)},\]
we are going to use global information of $Z$.  Following \cite[\S 4.2]{IM}, we have a map
\[
H^1(X,\Om^2_X)\to H^0(Z,\Om^1_Z),\]
which is obtained by taking the differential of a map
\[
\Def(X,Y)\to \Def(f,Y^{(2)})\]
from the space of deformations of $X$ keeping $Y$ fixed to the space of semi-universal deformations of $f$.
Using the symplectic form $\psi$, we have an isomorphism $N_{f(Z)/\Ypd; f(z)} \cong \Om^1_{Z,z}$.

From this, and the fact that the map $\beta$ is induced by deforming $X$ keeping $Y$ fixed,  we have a factorization
\[
\xymatrix{H^{2,1}(X)\ar[r]\ar[dr]_\beta& H^0(Z,\Om^1_Z)\ar@{->>}[d]\\
&N_{f(Z)/\Ypd;f(z)}}\]
where the vertical arrow is surjective.  Hence in this case $\beta$ is surjective. With what was said above this implies that for the Fano surface $d\AJ_X$ is an isomorphism.
\end{Examo}

It is a well-known standard fact that $H^1(\Om^2_X)\cong V$, the map being given for $P\in V$  by
\[
P\to \om_P = \Res_X \Om_P\]
where $\Om_P=\frac{P\Om}{F^2}$ with $\Om=\sum (-1)^ix_i dx_0\wedge\cdots \wedge \hat{dx}_i \wedge\cdots \wedge dx_4\in H^0(\Om^4_{\P^4}(5))$ being the standard form on $\P^4$.  If $z\in Z$ corresponds to the line $L_z\subset X$ where $L_z\in G(2,V^\ast)$, then for $\nu\in H^0(N_{L_z/X})\cong T_z Z$ it is direct to check that $$P\in L_z^\bot\subset V\implies \lra{d\AJ_X(\nu),P}=\int_{L_z}\nu\rfloor \om_P=0.$$  Since we have seen above that the induced map 
\[
d\AJ^\ast_X: V\to T^\ast_{L_z}Z\]
is surjective, it follows that $T_z Z\cong V^\ast /L_z$; i.e., 
\[
TZ=\lrc{ \bmp{2}{restriction to $Z\subset G(2,V^\ast)$ of the universal sub-bundle}}\]
which is the tangent bundle theorem for the Fano surface  (\cite{CG}).  In fact the above argument says that from $\ker\{H^0(\Om^1_Z)\to T^\ast_z Z\}$ we may recover $L_z$ and therefore the mapping
\[
J(X)\to \Alb(Z)\]
is an isomorphism.

\begin{rem}
Geometric arguments that $\alpha$ is injective and $\beta$ is surjective may be given as follows:
\beb
\item For $\alpha$ each line $L_z$ meets $Y$ in two points $L_z\cap Y$, and as $L_z$ varies one of these points will also vary (through a general point  of $X$ there are six lines);
\item For  a description of $\beta$, $\eta\in H^1(\Om^2_X)\cong V$ gives a first order deformation $X'$ of $X$ keeping  $Y$ fixed.  The line $L_z$ deforms (non-uniquely) to a line $L'$, and since $Y$ remains fixed
\[
L'\cap Y\in T_{f(z)}Y^{(2)}/\alpha(T_z(Z))=N_{f(Z)/Y^{(2)},f(z)}\]
is well defined.  We want to show that the $L'\cap Y$'s fill out the normal space.  Among the deformations of $X$ keeping fixed a  $Y$ given by $ F=0,Q=0 $ are the sections of $H^0(N_{X/\P V^\ast})\cong H^0(\cO_X(3))$ given by $P\cdot Q$ where $P\in H^0(\cO_X(1))=V$.  Since $P$ is free to vary this proves the surjectivity of the above map to the  normal space.
\eeb
\end{rem}
\begin{Examo}
Let $X_{2g-2}\subset \P^{g+1}$ be an anticanonically embedded Fano threefold of genus $g\geqq 5$, $-K_X$ non-divisible and Picard number one.  It is known that a general conic $C$ on $X$ has normal bundle isomorphic to $\cO_{\P^1}\oplus \cO_{\P^1}$; see \cite[Prop.\ 4.2.5]{IP}.  Then as before, we obtain that for a general conic $C$, the map $\alpha$ is injective.  If we let $g=6$, \cite{DIM} proves that for $X_{10}$, the variety of conics $Z$ contained in $X$ has a cotangent bundle that is generically globally generated.

Using the argument of Example 1, we obtain that for a general conic $C$ in $X_{10}$, $\beta$ is surjective, so that once again we may conclude that $d\AJ_X$ is an isomorphism.
\end{Examo}
\section{Discussion of the Proof of Theorem B}

 {\em Proof of \eqref{1.7}.} Referring to \eqref{4.1} we have
\[
\dim\fZ=h^0(N_{C/X})+\dim(\ker\delta).\]
From the commutative diagram
\[\begin{picture}(5,5) \put(80,-23){$\sidecong$}\end{picture}\xymatrix{
0\ar[r]& H^0(N_{C/X}\otimes K_X)\ar[r]& H^0(N_{C/X})\ar[d]^{d\AJ_X}\\
&H^1(N_{C/X})^\ast\ar[r]^{\delta^\ast}&H^{1,2}(X)}
\]
we have the identification \eqref{4.6}
\[
d\AJ_X\big|_{\ker \alpha}=\delta^\ast.\]
Thus the mapping 
\[
H^1(N_{C/X})^\ast\to T\Ypd \times T\cD_\la\]
has kernel equal to $\ker\delta^\ast$.  It follows that
\begin{align*}
\dim (F_\ast T\fZ)&=\dim T\fZ-\dim\ker\delta^\ast\\
&= d+h^1(N_{C/X})+\dim \ker\delta-\dim\ker\delta^\ast.\end{align*}
Using 
\[
\ker\delta^\ast = (\Rim \delta)^\bot\]
gives
\begin{align*}
\dim(\Rim \delta)^\bot &= h^1 (N_{C/X}) - \dim \Rim \delta\\
\dim (F_\ast T\fZ)&=d+\dim\ker\delta +\dim\Rim \delta= d+h^{2,1}(X)\\
&=\lrp{\frac12} \dim (\Ypd \times \cD_\la).\end{align*}

This argument shows that
 $d\AJ_X (\ker (H^0(N_{C/X}) \to T\Ypd) $  exactly compensates for the obstructions to freely deforming $C$ as $X$ deforms keeping $Y$ fixed.  
 
In summary, \emph{the symplectic geometry of $F:\fZ\to \Ypd \times \cD_\la$ exactly captures the infinitesimal Abel-Jacobi maps of the family of curves obtained by deforming the $C_z\subset X$  where $X$ deforms keeping $Y$ fixed.}

As mentioned, the complete proof of a more general version of Theorem B will be given elsewhere.

\bibliographystyle{amsalpha}

\bibliography{PG.bib}

\providecommand{\bysame}{\leavevmode\hbox to3em{\hrulefill}\thinspace}
\providecommand{\MR}{\relax\ifhmode\unskip\space\fi MR }
\providecommand{\MRhref}[2]{%
  \href{http://www.ams.org/mathscinet-getitem?mr=#1}{#2}
}
\providecommand{\href}[2]{#2}
\begin{thebibliography}{DIM12}

\bibitem[Bea04]{B}
Arnaud Beauville, \emph{Fano threefolds and {$K3$} surfaces}, The {F}ano
  {C}onference, Univ. Torino, Turin, 2004, pp.~175--184. \MR{2112574}

\bibitem[CG72]{CG}
C.~Herbert Clemens and Phillip~A. Griffiths, \emph{The intermediate {J}acobian
  of the cubic threefold}, Ann. of Math. (2) \textbf{95} (1972), 281--356.
  \MR{302652}

\bibitem[DIM12]{DIM}
O.~Debarre, A.~Iliev, and L.~Manivel, \emph{On the period map for prime {F}ano
  threefolds of degree 10}, J. Alg. Geom. \textbf{21} (2012), 21--59.

\bibitem[DM96]{DM}
Ron Donagi and Eyal Markman, \emph{Spectral covers, algebraically completely
  integrable, {H}amiltonian systems, and moduli of bundles}, Integrable systems
  and quantum groups ({M}ontecatini {T}erme, 1993), Lecture Notes in Math.,
  vol. 1620, Springer, Berlin, 1996, pp.~1--119. \MR{1397273}

\bibitem[IM07]{IM}
A.~Iliev and L.~Manivel, \emph{Prime {F}ano threefolds and integrable systems},
  Math. Ann. \textbf{339} (2007), 937--955.

\bibitem[IP99]{IP}
V.~Iskoviskivh and Y.~Prokhorov, \emph{Algebraic {G}eometry {V: F}ano
  {V}arieties}, Springer-Verlag, 1999.

\bibitem[Mar08]{M}
D.~Markushevich, \emph{An integrable system of {$K3$}-{F}ano flags}, Math. Ann.
  \textbf{342} (2008), no.~1, 145--156. \MR{2415319}

\end{thebibliography}

       \end{document}